\documentclass
[aps,showkeys,showpacs,notitlepage,superscriptaddress]{revtex4-1}
\pdfoutput=1
\usepackage{amsmath}
\usepackage{amssymb}
\usepackage{amsfonts}
\usepackage{graphicx}
\usepackage{hyperref}
\usepackage{color}
\usepackage{multirow}
\hypersetup{pdftitle=Comparative analysis of two discretizations of Ricci curvature for complex networks}

\usepackage{mathtools}

\begin{document}
\title{Comparative analysis of two discretizations of Ricci curvature for complex networks}

\author{Areejit Samal}
\thanks{AS and RPS contributed equally to this work}
\affiliation{The Institute of Mathematical Sciences (IMSc), Homi Bhabha National Institute (HBNI), Chennai 600113 India}
\author{R.P. Sreejith}
\thanks{AS and RPS contributed equally to this work}
\affiliation{The Institute of Mathematical Sciences (IMSc), Homi Bhabha National Institute (HBNI), Chennai 600113 India}
\author{Jiao Gu}
\affiliation{Jiangnan University, Wuxi 214122 P.R. China}
\author{Shiping Liu}
\affiliation{School of Mathematical Sciences, University of Science and Technology of China, Hefei 230026 P.R. China}
\author{Emil Saucan}
\email{Corresponding author: semil@ee.technion.ac.il}
\affiliation{Department of Applied Mathematics, ORT Braude College, Karmiel 2161002 Israel}
\affiliation{Departments of Mathematics and Electrical Engineering, Technion, Israel Institute of Technology, Haifa 3200003 Israel}
\author{J\"urgen Jost}
\email{Corresponding author: jost@mis.mpg.de}
\affiliation{Max Planck Institute for Mathematics in the Sciences, Leipzig 04103 Germany}
\affiliation{The Santa Fe Institute, Santa Fe, New Mexico 87501 USA}

\begin{abstract}
We have performed an empirical comparison of two distinct notions of discrete Ricci curvature for graphs or networks, namely, the Forman-Ricci curvature and Ollivier-Ricci curvature. Importantly, these two discretizations of the Ricci curvature were developed based on different properties of the classical smooth notion, and thus, the two notions shed light on different aspects of network structure and behavior. Nevertheless, our extensive computational analysis in a wide range of both model and real-world networks shows that the two discretizations of Ricci curvature are highly correlated in many networks. Moreover, we show that if one considers the augmented Forman-Ricci curvature which also accounts for the two-dimensional simplicial complexes arising in graphs, the observed correlation between the two discretizations is even higher, especially, in real networks. Besides the potential theoretical implications of these observations, the close relationship between the two discretizations has practical implications whereby Forman-Ricci curvature can be employed in place of Ollivier-Ricci curvature for faster computation in larger real-world networks whenever coarse analysis suffices.
\end{abstract}

\maketitle

\section{Introduction}

One of the central quantities associated to a Riemannian metric is the Ricci tensor. In Einstein's field equations, the energy-momentum tensor yields the Ricci tensor, and this determines the metric of space-time. In Riemannian geometry, the importance of the Ricci tensor came to the fore in particular through the work of Gromov \cite{Gromov1981}. The Ricci flow, introduced by Hamilton \cite{Hamilton1986}, culminated in the work of Perelman \cite{Perelman2002,Perelman2003} which solved the Poincar\`e and the more general Geometrization Conjecture for three-dimensional manifolds. On the other hand, there have been important developments extending the notion of Ricci curvature axiomatically to metric spaces more general than Riemannian manifolds \cite{Bakry2014,Lott2009,Sturm2006}. More precisely, one identifies metric properties on a Riemannian manifold that can be formulated in terms of local quantities such as growth of volumes of distance balls, transportation distances between balls, divergence of geodesics, and meeting probabilities of coupled random walks. On Riemannian manifolds such local quantities are implied by, or even equivalent to, Ricci curvature inequalities. Moreover when such metric properties are satisfied on some metric space, one says that the space satisfies the corresponding generalized Ricci curvature inequality. This research paradigm has been remarkably successful, and the geometry of metric spaces with such inequalities is currently a very active and fertile field of mathematical research (see for instance \cite{Bauer2017}). Of course, on Riemannian manifolds various such properties are equivalent to Ricci curvature inequalities and therefore also to each other. However, when passing to a discrete, metric setting, each approach captures different aspects of the classical Ricci curvature and thus, the various discretizations need no longer be equivalent. One such approach to Ricci curvature inequalities is Ollivier's \cite{Ollivier2007,Ollivier2009,Ollivier2010,Ollivier2013} construction on metric spaces.

There is also an older line of research \cite{Stone1976} that searches for the discretization of Ricci curvature on graphs and more general objects with a combinatorial structure. Here, one has exact quantities rather than only inequalities as in the aforementioned research. One elegant approach is by Chow and Luo \cite{Chow2003} based on circle packings which lent itself to many practical applications in graphics, medical imaging and communication networks \cite{Jin2007,Gu2013,Gao2014}. On the other hand, Ollivier's \cite{Ollivier2007,Ollivier2009,Ollivier2010,Ollivier2013} discretization has proven to be suitable for modelling complex networks as well as rendering interesting theoretic results with potential of future applications \cite{Lin2010,Lin2011,Bauer2012,Jost2014,Loisel2014,Ni2015,Sandhu2015a}. Yet another approach to discretization of Ricci curvature on polyhedral complexes, and more generally, $CW$ complexes is due to Forman \cite{Forman2003}. In recent work \cite{Sreejith2016,Sreejith2016directed,Sreejith2017,Weber2017,Saucan2018}, we have introduced the Forman's \cite{Forman2003} discretization to the realm of graphs and have systematically explored the Forman-Ricci curvature in complex networks. A crucial advantage of Forman-Ricci curvature is that, while it also captures important geometric properties of networks, it is far simpler to evaluate on large networks than Ollivier-Ricci curvature \cite{Sreejith2016,Saucan2018}. In this contribution, we have performed an extensive empirical comparison of the Forman-Ricci curvature and Ollivier-Ricci curvature in complex networks. In addition, we have also performed an empirical analysis in complex networks of the augmented Forman-Ricci curvature which accounts for two-dimensional simplicial complexes arising in graphs. We find that the Forman-Ricci curvature, especially the augmented version, is highly correlated to Ollivier-Ricci curvature in many model and real networks. This renders Forman-Ricci curvature a preferential tool for the analysis of very large networks with various practical applications.

Although, in this contribution, we show that Forman-Ricci curvature is highly correlated to Ollivier-Ricci curvature in many networks, one should not construe from this observation that we introduce Forman-Ricci curvature as a substitute (and certainly not as a ``proxy'' \cite{Pal2018}) for Ollivier-Ricci curvature. As mentioned above, and as we shall further explain in the following section, the two discretizations of Ricci curvature capture quite different aspects of network behavior. Indeed the specific definitions of both Ollivier's and Forman's discretizations of Ricci curvature prescribe some of their respective essential properties that have important consequences in certain significant applications. Therefore, we shall detail these definitions and not restrict ourselves to the mere technical defining formulas.

Given that networks permeate almost every field of research \cite{Wasserman1994,Watts1998,Barabasi1999,Albert2002,Feng2007,Newman2010,Fortunato2010}, an important challenge has been to unravel the architecture of complex networks. In particular, the development of geometric tools \cite{Eckmann2002,Ollivier2009,Lin2010,Lin2011,Bauer2012,Jost2014,Wu2015,Ni2015,Sandhu2015a,Sreejith2016,Bianconi2017}, and mainly curvature, allow us to gain deep insights into the structure, dynamics and evolution of networks. It is in the very nature of discretization of differential geometric properties that each such discrete notion sheds a different light and understanding upon the studied object, for example, a network. In particular, Ollivier's curvature is related to clustering and network coherence via the distribution of the eigenvalues of the graph Laplacian, giving insights into the global and local structure of networks. In contrast, Forman's curvature captures the geodesics dispersal property and also gives information on the algebraic topological structure of the network. Furthermore, Forman's curvature is simple to compute and can easily be extended to analyze both directed networks and hyper-networks \cite{Sreejith2016,Sreejith2016directed,Sreejith2017,Weber2017}. Given the contrast between the two discretizations of Ricci curvature at hand, the empirically observed correlation in many networks is quite surprising and encouraging. Moreover, both types of curvature admit natural Ricci curvature flows \cite{Gu2013,Weber2017} that enable the study of long time evolution and prediction of networks. Moreover, the observed correlation further increases the relevance and importance of future investigation of discrete Ricci flows for the better understanding of the structure and evolution of complex networks.

Note that in Riemannian geometry, the Ricci tensor encodes all the essential properties of a Riemannian metric. Similarly, it is an emerging principle that Ricci curvature, because it evaluates edges instead of vertices, also captures the basic structural aspects of a network. Both Ollivier-Ricci curvature and Forman-Ricci curvature are {\em edge-based} measures which assign a number to each edge of a (possibly weighted and directed) network that encodes local geometric properties in the vicinity of that edge. We highlight that {\em edges} are what networks are made of as the {\em vertices} alone do not yet constitute a network.

\begin{figure}
\includegraphics[width=.7\columnwidth]{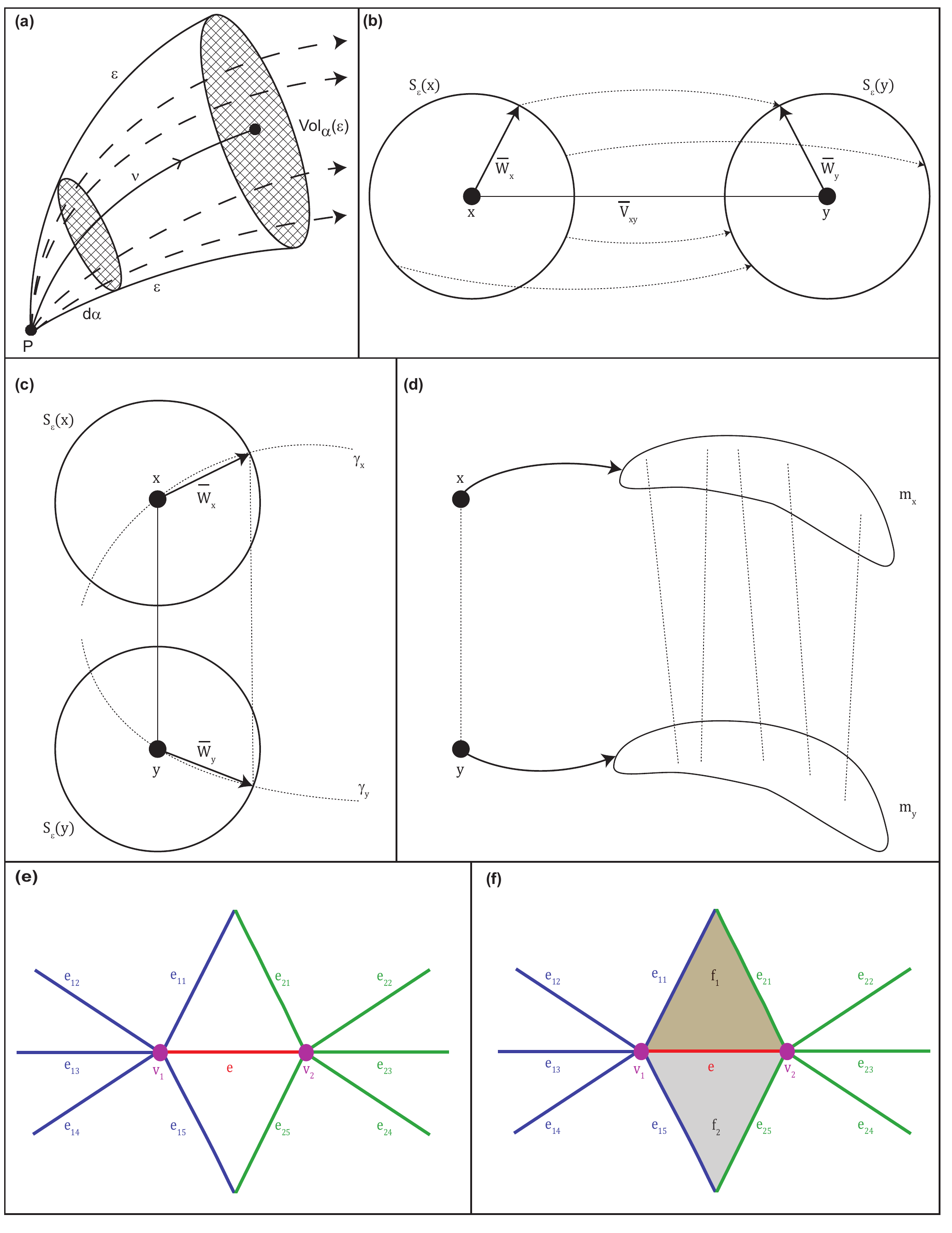}
\caption{{\bf (a)} The geometric interpretation of Ricci curvature. Ricci curvature measures the growth of volumes, more precisely, the growth of $(n-1)$-dimensional solid angles in the direction of the vector ${\bf v}$. It also measures the dispersion rate of the family of geodesics with the same initial point, that are contained within the given solid angle. {\bf (b)-(d)} The interpretation of Ollivier-Ricci curvature. {\bf (b)} Given two close points $x$ and $y$ in a Riemannian manifold of dimension $n$, defining a tangent vector $\bar{v}_{xy}$, one can consider the parallel transport in the direction $\bar{v}_{xy}$. Then points on a infinitesimal sphere $S_\varepsilon(x)$ centered at $x$, are transported to points on the corresponding sphere $S_\varepsilon(y)$ by a distance equal to $d(x,y)\left(1 - \frac{\varepsilon^2}{2n}\rm{Ric}(\bar{v}_{xy})\right)$, on the average. {\bf (c)} In Riemannian manifolds of positive (respectively, negative) curvature, balls are closer (respectively, farther) than their centers. Thus, in spaces of positive Ricci curvature spheres are closer than their centers, while in spaces of  negative curvature they are farther away. {\bf (d)} To generalize this idea to metric measure spaces, one has to replace the (volumes of) spheres or balls, by measures $m_x$, $m_y$. Points will be transported by a distance equal to $(1 - \kappa)d(x,y)$, on the average, where $\kappa = \kappa(x,y)$ represents the coarse (Ollivier) curvature along the geodesic segment $xy$. This illustration is an adaptation of the original figure \cite{Ollivier2013}. {\bf (e)-(f)} Forman-Ricci curvature of an edge $e$ connecting the vertices $v_1$ and $v_2$ and contributions from edges parallel to the edge $e$ under consideration. An edge is said to be parallel to a given edge $e$, if it has in common with $e$ either a \textit{child} (i.e., a lower dimensional face), or a \textit{parent} (i.e., a higher dimensional face), but not both simultaneously. In part {\bf (e)}, all the edges $e_{11}, \ldots, e_{15}$ are parallel to $e$ because they share the vertex $v_1$, while the edges $e_{21}, \ldots, e_{25}$ are parallel to $e$ because they share the vertex $v_2$. In contrast, in part {\bf (f)}, edges $e_{11}, e_{21}, e_{15}, e_{25}$ are not parallel anymore to the edge $e$, because they have common \textit{children} with $e$ (namely, $v_1$ and $v_2$) and a common parent with $e$ (namely, $f_1$ or $f_2$). In consequence, edges $e_{11}, e_{21}, e_{15}, e_{25}$ do not contribute in the computation of the Augmented Forman-Ricci curvature of edge $e$ which also accounts for the two-dimensional simplicial complexes $f_1$ and $f_2$.}
\label{fig:ricciinterpret}
\end{figure}

\section{Discrete Ricci curvatures on networks}

We briefly present here the geometric meaning of the notion of Ricci curvature, as well as the two discretizations considered herein. For other discretizations of this type of curvature and their applications, see for instance \cite{Gu2013}.

\subsection{Ricci curvature}

In Riemannian geometry curvature measures the deviation of the manifold from being locally Euclidean. Ricci curvature quantifies that deviation for tangent directions. It controls the average dispersion of geodesics around that direction. It also controls the growth of the volume of distance balls and spheres. In fact, these two properties are related, as can be seen from the following formula \cite{Heintze1978}:	
\begin{equation}
\label{eq:Ricci}
{\rm Vol}_\alpha(\varepsilon) = d\alpha\,\varepsilon^{n-1}\left(1 - \frac{{\rm Ric}({\bf v})}{3}\varepsilon^2 + o(\varepsilon^2)\right)\,.
\end{equation}
Here, $n$ is the dimension of the Riemannian manifold in question, and ${\rm Vol}_\alpha(\varepsilon)$ is the $(n-1)$-volume generated within an $n$-solid angle $d\alpha$ by geodesics of length $\varepsilon$ in the direction of the vector ${\bf v}$ (i.e., it controls the growth of measured angles). Thus, Ricci curvature controls both divergence of geodesics and volume growth (Figure \ref{fig:ricciinterpret}(a)). In dimension $n=2$, Ricci curvature reduces to the classical {\it Gauss curvature}, and can therefore be easily visualized.

As we shall see, the two discretizations of Ricci curvature by Ollivier and Forman considered here for networks capture different properties of the classical (smooth) notion. Forman's definition expresses dispersal (diffusion), while Ollivier's definition compares the averaged distance between balls to the distance between their centers. Thus, the two definitions lead to different generalization of classical results regarding Ricci curvature. In this respect, Ollivier's version seems to be advantageous, since, in addition to certain geometric properties, analytic inequalities also hold, whereas Forman's version encapsulates mainly the topology of the underlying space.

Nevertheless, in our specific context of complex networks, as we shall show in the sequel, the definitions by Ollivier and Forman are highly correlated in many networks. Therefore, for the empirical analysis of large networks, at least in a first approximation, from the analysis of Forman's definition, one can also make inferences about the properties encoded by the Ollivier's definition. For instance, Ollivier's curvature is, by its very definition, excellently suited to capture diffusion and stochastic properties of a given network. Unfortunately, the computation of Ollivier-Ricci curvature might be prohibitive for many large complex networks. In contrast, due to its simple, combinatorial formula, Forman-Ricci curvature is easy and fast to compute \cite{Sreejith2016}. Given the basic equivalence, at least on a statistical level, between these two discretizations, one can therefore determine, at least in first approximation, many properties encapsulated by Ollivier's curvature via simple computations with Forman's curvature. However, for a finer analysis, each of the two discrete Ricci curvatures should be employed in the context that best befits the geometrical phenomenology it encapsulates.


\subsection{Ollivier-Ricci curvature}

Ollivier's approach \cite{Ollivier2007,Ollivier2009,Ollivier2010,Ollivier2013} interprets Eq. \ref{eq:Ricci} as follows: If a small ball $B_x$ of radius $\varepsilon$ and centered at $x$ is mapped, via parallel transport \cite{Jost2017} to a corresponding ball $B_y$ centered at $y$, then the average distance between points on $B_x$ and their corresponding points on $B_y$ is:
\begin{equation}
\label{eq:Ollivier1}
\delta\left(1 - \frac{\varepsilon^2}{2(n+2)}{\rm Ric}(v) + O(\varepsilon^3 + \varepsilon^2\delta) \right)\,,
\end{equation}
where $d(x,y) = \delta$, and where $\varepsilon, \delta \rightarrow 0$. Thus, we can {\em synthetically} characterize Ollivier-Ricci curvature \cite{Ollivier2013} by the following phrase: ``In positive (negative) curvature, balls are closer (farther) than their centers are''. Balls are given by their volume measures, and in fact, one may define a transportation distance for any two (normalized) measures. In this sense, Ollivier's notion compares the distance between the centers of their balls with that between their measures (Figure \ref{fig:ricciinterpret}(b)-(d)). For the distance between the centers one takes (of course) the given metric of the underlying space, i.e., manifold, mesh, network, etc. As for the distance between measures, there is a natural choice, the Wasserstein transportation metric $W_1$ \cite{Vaserstein1969}. More formally, Ollivier's curvature is defined as:
\begin{equation}
\label{eq:Ollivier2}
\kappa(x,y) = 1 - \frac{W_1(m_x,m_y)}{d(x,y)}\,;
\end{equation}
where $m_x, m_y$ represent the measures of the balls around $x$ and $y$, respectively. Here, given that the measure $m$, associated to the discrete set of vertices of a graph (network) is obviously a discrete measure, the Wasserstein distance $W_1(m_x, m_y)$, i.e. the transportation distance between the two probability measures $m_x$ and $m_y$, is given by
\begin{equation}
\label{eq:Ollivier3}
W_1(m_x, m_y)=\inf_{\mu_{x,y}\in \prod(m_x, m_y)}\sum_{(x',y')\in V\times V}d(x', y')\mu_{x,y}(x', y'),
\end{equation}
with $\prod(m_x, m_y)$ being the set of probability measures $\mu_{x,y}$ that satisfy:
\begin{equation}
\label{eq:Ollivier4}
\sum_{y'\in V}\mu_{x,y}(x', y')=m_x(x'), \,\,\sum_{x'\in V}\mu_{x,y}(x', y')=m_y(y').
\end{equation}
Measures satisfying Eq. \ref{eq:Ollivier4} start with the measure $m_x$ and end up with $m_y$, and represent all the transportation possibilities of the mass (measure) $m_x$ to the measure $m_y$, by \textit{disassembling} it, transporting it, along all possible paths, and \textit{reassembling} it as $m_y$. $W_1(m_x, m_y)$ is the minimal cost (measured in terms of distances) to transport the mass of $m_x$ to that of $m_y$. Note that the distance $d$ in Eq. \ref{eq:Ollivier3} above can be any useful or expressive graph metric. However, in practice, when considering the Wasserstein metric and Ollivier-Ricci curvature for unweighted networks, the {\it combinatorial} metric is naturally considered.

In the Riemannian setting, Ollivier's definition  reduces to the classical one. More precisely, if $M^n$ is a Riemannian manifold, with its natural measure $d{\rm Vol}$, then for $d(x,y)$ small enough and $v$ the unit tangent vector at $x$ on the geodesic $\overline{xy}$
\begin{equation}
\kappa(x,y) =  \frac{\varepsilon^2}{2(n+2)}{\rm Ric}(v) +  O(\varepsilon^3 + \varepsilon^2d(x,y))\,.
\end{equation}

The Wasserstein distance \cite{Vaserstein1969} between two vertices in a network depends on the triangles, quadrangles and pentagons that they are contained in (see for instance \cite{Jost2014,Bhattacharya2015}). It can also be computed in terms of random walks on a graph, where one has the choice between the {\em lazy } \cite{Ni2015} and the {\em non-lazy} \cite{Jost2014} random walk. While the two variants are clearly equivalent from a theoretical viewpoint, the choices may render differences in the implementation. In this work, we have used the lazy random walk option within the open-source implementation of Ollivier-Ricci curvature, originally developed by P. Romon and improved by E. Madsen, within SageMath software (http://www.sagemath.org/) for our computations.

While Ollivier-Ricci curvature is essentially defined on edges, one can define Ollivier-Ricci curvature of a vertex \cite{Sandhu2015a} as the sum of the Ollivier-Ricci curvatures of edges incident on that vertex in the network, and this is analogous to scalar curvature in Riemannian geometry \cite{Jost2017}.


\subsection{Forman-Ricci curvature}

Forman's definition is conceptually quite different from Ollivier's definition. To begin with, Forman's definition works in the framework of \textit{weighted $CW$ cell complexes}, rather than that of Markov chains and metric measure spaces, as Ollivier's definition does. The weighted $CW$ cell complexes are of fundamental importance in topology and include both polygonal meshes and weighted graphs. In the setting of weighted $CW$ cell complexes, Forman's definition develops an abstract version of a classical formula in differential geometry or geometric analysis, the so called \textit{Bochner-Weitzenb\"{o}ck formula} (see for instance \cite{Jost2017}), that relates curvature to the classical (Riemannian) Laplace operator.

Forman \cite{Forman2003} derived an analogue of the Bochner-Weitzenb\"{o}ck formula that holds in the setting of $CW$ complexes. In the 1-dimensional case, i.e. of graphs or networks, it takes the following form \cite{Sreejith2016}:
\begin{equation}
\label{eq:Forman1}
{\rm F}(e) = w_e \left( \frac{w_{v_1}}{w_e} +  \frac{w_{v_2}}{w_e}  - \sum_{e_{v_1}\ \sim\ e,\ e_{v_2}\ \sim\ e} \left[\frac{w_{v_1}}{\sqrt{w_e w_{e_{v_1} }}} + \frac{w_{v_2}}{\sqrt{w_e w_{e_{v_2} }}} \right] \right)\,
\end{equation}
where $e$ denotes the edge under consideration between two nodes $v_1$ and $v_2$, $w_e$ denotes the weight of the edge $e$ under consideration, $w_{v_1}$ and $w_{v_2}$ denote the weights associated with the vertices $v_1$ and $v_2$, respectively, $e_{v_1} \sim e$ and $e_{v_2} \sim e$ denote the set of edges incident on vertices $v_1$ and $v_2$, respectively, after excluding the edge $e$ under consideration which connects the two vertices $v_1$ and $v_2$ (Figure \ref{fig:ricciinterpret}(e)). Since edges in the discrete setting of networks naturally correspond to vectors or directions in the smooth context, the above formula represents, in view of the classical Bochner-Weitzenb\"{o}ck formula, a discretization of Ricci curvature. For gaining further intuition regarding this discretization of Ricci curvature in its generality, the reader is referred to Forman's original work \cite{Forman2003}, and to our previous papers \cite{Sreejith2016,Saucan2018} for more insight on its adaptation to networks.

In the combinatorial case, i.e. for $w_e = w_v = 1, \; e \in E(G), v \in V(G)$, where $E(G)$ and $V(G)$ represent the set of edges and vertices, respectively, in graph $G$, the above formula (Eq. \ref{eq:Forman1}) reduces to the quite simple and intuitive expression:
\begin{equation}
\label{eq:Forman2}
{\rm F} (e) = 4 - \sum_{v \sim e} \deg(v) \;
\end{equation}
where $v \sim e$ denote the vertices anchoring the edge $e$. This simple case captures the role of Ricci curvature as a measure of the flow through an edge and illustrates how Ricci curvature captures the \textit{social behavior} of geodesics dispersal depicted in Figure \ref{fig:ricciinterpret}. While Forman-Ricci curvature is essentially defined on edges, one can easily define Forman-Ricci curvature of a vertex \cite{Sreejith2017} as the sum of the Forman-Ricci curvatures of edges incident on that vertex in the network.

\subsubsection*{Augmented Forman-Ricci curvature}

From a graph, one may construct two-dimensional polyhedral complexes by inserting a two-dimensional simplex into any connected triple of vertices (or cycle of length 3), a tetragon into any cycle of length 4, a pentagon into a cycle of length 5, and so on. This is natural, if, for instance, one wants to represent higher order correlations between vertices in the network. Again, Forman's scheme assigns a Ricci curvature to such a complex, via the following formula, which also includes possible weights $w$ of simplices, edges, and vertices:
\begin{equation}
\label{eq:Forman3}
{\rm F}^{\#} (e) = w_e \left[ \left( \sum_{e < f} \frac{w_e}{w_f}+\sum_{v < e} \frac{w_v}{w_e} \right) \right. - \left. \sum_{\hat{e} \parallel e} \left| \sum_{\hat{e},e < f} \frac{\sqrt{w_e \cdot w_{\hat{e}}}}{w_f} - \sum_{v 	 < e, v < \hat{e}} \frac{w_v}{\sqrt{w_e \cdot w_{\hat{e}}}} \right| \right] \; ;
\end{equation}
where $w_e$ denotes weight of edge $e$, $w_v$ denotes weight of vertex $v$, $w_f$ denotes weight of face $f$, $\sigma < \tau$ means that $\sigma$ is a face of $\tau$, and where $||$ signifies \textit{parallelism}, i.e. the two cells have a common \textit{parent} (higher dimensional face) or a common \textit{child} (lower dimensional face), but not both a common parent and common child. In particular, we have employed Eq. \ref{eq:Forman3} to define an \textit{Augmented Forman-Ricci curvature} of an edge which also accounts for two-dimensional simplicial complexes or cycles of length 3 arising in graphs while neglecting cycles of length 4 and greater (Figure \ref{fig:ricciinterpret}(f)).

In unweighted networks, $w_f = w_e = w_v = 1, \; \forall f \in F(G), e \in E(G), v \in V(G)$, where $F(G)$, $E(G)$ and $V(G)$ represent the set of faces, edges and vertices, respectively, in graph $G$. In such unweighted networks, we remark that there is a simple relationship \cite{Webere2017triangle} between Forman-Ricci curvature ${\mathrm F} (e)$ of an edge $e$ and Augmented Forman-Ricci curvature ${\mathrm F}^{\#} (e)$  of an edge $e$, namely,
\begin{equation}
\label{eq:Forman4}
{\rm F}^{\#} (e) = {\rm F} (e) + 3m
\end{equation}
where $m$ is the number of triangles containing edge $e$ under consideration in the network. In this work, we have explored both Forman-Ricci curvature and its augmented version in model and real-world networks.


\subsection{Ollivier's vs. Forman's Ricci curvature: A first comparison}

As we have seen in detail in the previous section, and already explained in the Introduction, the two types of discrete Ricci curvature, Ollivier's and Forman's, express different geometric properties of a network, and they can therefore be quite different from each other for specific networks. In this section, let us consider some simple examples.

As the first example, consider a complete graph on $n$ vertices. Then any two vertices share $n-2$ neighbors in the complete graph, and therefore, the corresponding balls largely overlap. The transportation distance between the balls is thus very small in a complete graph, and thus, the Ollivier-Ricci curvature (Eq. \ref{eq:Ollivier2}) is almost 1 for large $n$, the largest possible value. On the other hand, the degree of any vertex is $n-1$ in a complete graph, and therefore, the Forman-Ricci curvature (Eq. \ref{eq:Forman2}) takes the most negative possible value. Thus, for such complete graphs, the two types of Ricci curvature behave in opposite fashion. The reason is that Ollivier-Ricci curvature is positively affected by triangles whereas Forman-Ricci curvature is not at all. Thus, it is not surprising that locally they can numerically diverge from each other. As the second example, consider a star graph, that is, a graph consisting of a central vertex $v_0$ that is connected to all other vertices $v_1,\dots ,v_m$, while these vertices have no further connections. Consider an edge, for example, $e=(v_0,v_1)$ in the star graph. The neighborhood of $v_1$ consists of $v_0$ only, while that of $v_0$ contains all the vertices $v_1,\dots ,v_m$ in the star graph. Since each of these vertices $v_1,\dots ,v_m$ have distance $1$ from $v_0$ in the star graph, the transportation cost is $1$, and hence the Ollivier-Ricci curvature is $0$. In this example of a star graph, there are no triangles. In contrast, the Forman-Ricci curvature of the edge in the star graph is $3-m$. As the third example, consider a double star graph, that is, take two stars with vertices $v_0, v_1,\dots ,v_m$ and $v'_0, v'_1,\dots ,v'_{m'}$, where the two central vertices $v_0$ and $v'_0$ of the stars are connected by an edge. In this case of double star graph, almost all vertices in their respective neighborhoods are a distance $3$ apart, and so, the Ollivier-Ricci curvature of the edge $(v_0,v'_0)$ is quite negative, and so is the Forman-Ricci curvature, which equals $2-m -m'$. Thus, the second example of a star graph is an intermediate between the first example of a complete graph and the third example of a double star graph.

While these examples suggest an equivocal picture wherein sometimes the two discretizations of Ricci curvatures are aligned, but in other cases, they may show an opposite behavior, our numerical results in complex networks which are reported in the following sections show that, Ollivier-Ricci and Forman-Ricci curvature in many networks are highly correlated to each other. Thus, in several model and real networks that we have investigated, large degrees of the vertices bounding an edge do not correlate highly with large fractions of triangles or other short loops containing these vertices. Furthermore, if we augment the definition of the Forman-Ricci curvature to account for two-dimensional simplicial complexes (i.e., triads or cycles of length 3) arising in graphs (Eqs. \ref{eq:Forman3} and \ref{eq:Forman4}), then such an Augmented Forman-Ricci curvature is even better correlated at \textit{small scale} to Ollivier-Ricci curvature, as in the augmented definition the triangles of vertices no longer contribute negatively to Forman-Ricci curvature. In the sequel, we shall also show that the Augmented Forman-Ricci curvature is better correlated to Ollivier-Ricci curvature in both model and real-world networks.


\section{Benchmark dataset of complex networks}
\label{dataset}

We have considered four models of undirected networks, namely, Erd\"{o}s-R\'{e}nyi (ER) \cite{Erdos1961}, Watts-Strogatz (WS) \cite{Watts1998}, Barab\'{a}si-Albert (BA) \cite{Barabasi1999} and Hyperbolic Graph Generator (HGG) \cite{Krioukov2010}. The ER model \cite{Erdos1961} produces an ensemble of random graphs $G(n,p)$ where $n$ is the number of vertices and $p$ is the probability that each possible edge exists between any pair of vertices in the network. The WS model \cite{Watts1998} generates small-world networks which exhibit both a high clustering coefficient and a small average path length. In the WS model, an initial regular graph is generated with $n$ vertices on a ring lattice with each vertex connected to its $k$ nearest neighbours. Subsequently an endpoint of each edge in the regular ring graph is rewired with probability $\beta$ to a new vertex selected from all the vertices in the network with a uniform probability. The BA model \cite{Barabasi1999} generates scale-free networks which exhibit a power-law degree distribution. In the BA model, an initial graph is generated with $m_0$ vertices. Thereafter, a new vertex is added to the initial graph at each step of this evolving network model such that the new vertex is connected to $m$ $\le$ $m_0$ existing vertices, selected  with a probability proportional to their degree. Thus, the BA model implements a preferential attachment scheme whereby high-degree vertices have a higher chance of acquiring new edges than low-degree vertices. The HGG model \cite{Krioukov2010,Aldecoa2015} can produce random hyperbolic graphs with power-law degree distribution and non-vanishing clustering. In the HGG model, the $n$ vertices of the network are placed randomly on a hyperbolic disk, and thereafter, pairs of vertices are connected based on some probability which depends on the hyperbolic distance between vertices. In the HGG model, the input parameters \cite{Krioukov2010,Aldecoa2015} are the number of vertices $n$, the target average degree $k$, the target exponent $\gamma$ of the power-law degree distribution and temperature $T$. In this work, we have used HGG model with default input parameters of $\gamma=2$  and $T=0$ to generate hyperbolic random geometric graphs. Note that the input parameters, $\gamma$ and $T$, of the HGG model \cite{Krioukov2010,Aldecoa2015} can be varied to produce other random graph ensembles such as configuration model, random geometric graphs on a circle and ER graphs.

Supplementary Table S1 lists the model networks analyzed in this work along with the number of vertices, number of edges, average degree and edge density of each network. In each model, we have chosen different combinations of input parameters to generate networks with different sizes and average degree (Supplementary Table S1). Moreover, we have sampled 100 networks starting with different random seed for a specific combination of input parameters from each generative model, and the results reported in the next section for model networks in an average over the sample of 100 networks with chosen input parameters (Supplementary Tables S2-S5).

We have also considered seventeen widely-studied real undirected networks. These are six communication or infrastructure networks, the Chicago road network \cite{Eash1979}, the Euro road network \cite{Subelj2011}, the US Power Grid network \cite{Leskovec2007}, the Contiguous US States network \cite{Knuth2005}, the autonomous systems network \cite{Leskovec2007} and an Email communication network \cite{Guimera2003}. In the Chicago road network, the 1467 vertices correspond to transportation zones within the Chicago region, and the 1298 edges are roads in the region linking them. In the Euro road network, the 1174 vertices are cities in Europe, and the 1417 edges are roads in the international E-road network linking them. In the US Power Grid network, the 4941 vertices are generators or transformers or substations in the western states of the USA, and the 6594 edges are power supply lines linking them. In the Contiguous US States network, the 48 vertices correspond to the 48 contiguous states of USA (except the two states, Alaska and Hawaii, which are not connected by land with the other 48 states), and the 107 edges represent land border between the states. In the autonomous systems (AS) network, the 26475 vertices are autonomous systems of the Internet, and the 53381 edges represent communication between autonomous systems connected to each other from the CAIDA project. In the Email communication network, the 1133 vertices are users in the University Rovira i Virgili in Tarragona in Spain, and the 5451 edges represent direct communication between them. We have considered five social networks, the Zachary karate club \cite{Zachary1977}, the Jazz musicians network \cite{Gleiser2003}, the Hamsterster friendship network, the Dolphin network \cite{Lusseau2003} and the Zebra network \cite{Sundaresan2007}. In the Zachary karate club, the 34 vertices correspond to members of an university karate club, and the 78 edges represent ties between members of the club. In the Jazz musicians network, the 198 vertices correspond to Jazz musicians, and the 2742 edges represent collaboration between musicians. In the Hamsterster friendship network, the 2426 vertices are users of hamsterster.com, and the 16631 edges represent friendship or family links between them. In the Dolphin network, the 62 vertices correspond to bottlenose Dolphins living off Doubtful Sound in South West New Zealand, and the 159 edges represent frequent associations among Dolphins observed between 1994 and 2001. In the Zebra network, the 27 vertices correspond to Grevy's Zebras in Kenya, and the 111 edges represent observed interaction between Zebras during the study \cite{Sundaresan2007}. We have also considered a scientific co-authorship network based on papers from the arXiv's Astrophysics (astro-ph) section \cite{Leskovec2007} where the 18771 vertices correspond to authors and the 198050 edges represent common publications among authors. We have also considered the PGP network \cite{Boguna2004}, an online contact network, where the 10680 vertices are users of the Pretty Good Privacy (PGP) algorithm, and the 24316 edges represent interactions between the users. We have also considered a linguistic network, an adjective-noun adjacency network \cite{Newman2006}, where the 112 vertices are nouns or adjectives, and the 425 edges represent their presence in adjacent positions in the novel David Copperfield by Charles Dickens. We have considered three biological networks, the yeast protein interaction network \cite{Jeong2001}, the PDZ domain interaction network \cite{Beuming2005} and the human protein interaction network \cite{Rual2005}. In the yeast protein interaction network, the 1870 vertices are proteins in yeast \textit{Saccharomyces cerevisiae}, and the 2277 edges are interactions between them. In the PDZ domain interaction network, the 212 vertices are proteins, and the 244 edges are PDZ-domain mediated interactions between proteins. In the human protein interaction network, the 3133 vertices are proteins, and the 6726 edges are interactions between human proteins as captured in an earlier release of the proteome-scale map of human binary protein interactions. The seventeen empirical networks analyzed here were downloaded from the KONECT \cite{Kunegis2013} database. Supplementary Table S1 lists the real networks analyzed in this work along with number of vertices, number of edges, average degree and edge density of each network.

We remark that the above-mentioned model and real-world networks considered in this work are unweighted graphs, and thus, the weights of vertices, edges and two-dimensional simplicial complexes are taken to be 1 while computing the Forman-Ricci curvature and its augmented version. Furthermore, the largest connected component of the above-mentioned model and real-world networks is considered while computing the Ollivier-Ricci curvature of edges. In earlier work \cite{Sreejith2016,Sreejith2017}, we had characterized the Forman-Ricci curvature of edges and vertices in some of the above-mentioned networks. In the present work, we have compared the Forman-Ricci curvature and its augmented version with Ollivier-Ricci curvature in above-mentioned networks.


\section{Results and Discussion}
\label{results}

\subsection{Comparison between Forman-Ricci and Ollivier-Ricci curvature in model and real networks}

We have compared the Ollivier-Ricci with Forman-Ricci and Augmented Forman-Ricci curvature of edges in model networks (Table \ref{tab:ORFRedge} and Supplementary Table S2). In random ER networks, small-world WS networks and scale-free BA networks, we find a high positive correlation between the Ollivier-Ricci and Forman-Ricci curvature of edges or between Ollivier-Ricci and Augmented Forman-Ricci curvature of edges when the model networks are sparse with small average degree, however, the observed correlation vanishes with increase in average degree of model networks (Table \ref{tab:ORFRedge} and Supplementary Table S2). In hyperbolic random geometric graphs, we also find a high positive correlation between the Ollivier-Ricci and Forman-Ricci curvature of edges or between Ollivier-Ricci and Augmented Forman-Ricci curvature of edges, however, the observed correlation in the hyperbolic graphs seems relatively less dependent on average degree of networks based on our limited exploration of the parameter space (Table \ref{tab:ORFRedge} and Supplementary Table S2). We remark that hyperbolic random geometric graphs unlike ER, WS and BA networks have explicit geometric structure. Note that the Augmented Forman-Ricci in comparison to Forman-Ricci curvature of edges has typically higher positive correlation with Ollivier-Ricci curvature of edges in ER, WS and BA models (Table \ref{tab:ORFRedge} and Supplementary Table S2). Moreover, WS networks have higher clustering coefficient (and thus, higher proportion of triads) in comparison to ER or BA networks with same number of vertices and average degree, and thus, it is not surprising to observe that the Augmented Forman-Ricci curvature in comparison to Forman-Ricci curvature of edges has much higher positive correlation with Ollivier-Ricci curvature of edges in WS networks, especially, when networks become denser with increase in average degree (Table \ref{tab:ORFRedge} and Supplementary Table S2). This last result is expected because the Augmented Forman-Ricci curvature of edges also accounts for two-dimensional simplicial complexes or cycles of length 3 arising in graphs (see discussion in Theory section and Figure \ref{fig:ricciinterpret}(e)-(f)).

We have also compared the Ollivier-Ricci with Forman-Ricci and Augmented Forman-Ricci curvature of edges in seventeen real-world networks. In several of the analyzed real-world networks, we find a moderate to high positive correlation between Ollivier-Ricci and Forman-Ricci curvature of edges (Table \ref{tab:ORFRedge} and Supplementary Table S2). We highlight that some of the real-world networks such as Astrophysics co-authorship network, Email communication network, Jazz musicians network and Zebra network have very weak or no correlation between Ollivier-Ricci and Forman-Ricci curvature of edges (Table \ref{tab:ORFRedge} and Supplementary Table S2). However, in most real-world networks analyzed here, we find a moderate to high positive correlation between Augmented Forman-Ricci and Ollivier-Ricci curvature of edges (Table \ref{tab:ORFRedge} and Supplementary Table S2). Interestingly, we also find that the Augmented Forman-Ricci curvature has moderate to high correlation with Ollivier-Ricci curvature of edges in  Astrophysics co-authorship network, Email communication network, Jazz musicians network and Zebra network where Forman-Ricci curvature has very weak or no correlation with Ollivier-Ricci curvature of edges (Table \ref{tab:ORFRedge} and Supplementary Table S2). Thus, at the level of edges, we observe a positive correlation between Ollivier-Ricci and Forman-Ricci curvature, especially, the augmented version, in many networks (Table \ref{tab:ORFRedge} and Supplementary Table S2).

From the definition of the Ollivier-Ricci and Forman-Ricci curvature of edges, it is straightforward to define Ollivier-Ricci and Forman-Ricci curvature of vertices in networks \cite{Sandhu2015a,Sreejith2017} as the sum of the Ricci curvatures of the edges incident on the vertex in the network. Note that the definition of Ollivier-Ricci and Forman-Ricci curvature of vertices in networks \cite{Sandhu2015a,Sreejith2017} is a direct discrete analogue of the scalar curvature in Riemannian geometry \cite{Jost2017}.

We have compared the Ollivier-Ricci with Forman-Ricci and Augmented Forman-Ricci curvature of vertices in model networks (Table \ref{tab:ORFRvertex} and Supplementary Table S3).  In random ER networks, small-world WS networks and scale-free BA networks, we find a high positive correlation between the Ollivier-Ricci and Forman-Ricci curvature of vertices or between Ollivier-Ricci and Augmented Forman-Ricci curvature of vertices, and the observed correlation seems to have minor dependence on size or average degree of networks based on our limited exploration of the parameter space  (Table \ref{tab:ORFRvertex} and Supplementary Table S3). In most hyperbolic random geometric graphs analyzed here, we also find a moderate positive correlation between the Ollivier-Ricci and Forman-Ricci curvature of vertices or between Ollivier-Ricci and Augmented Forman-Ricci curvature of vertices (Table \ref{tab:ORFRvertex} and Supplementary Table S3). Note that in random ER networks, small-world WS networks and scale-free BA networks, the Spearman correlation is typically higher than Pearson correlation between Ollivier-Ricci and Forman-Ricci curvature of vertices, however, in the hyperbolic random geometric graphs, the Spearman correlation is typically lower than Pearson correlation between Ollivier-Ricci and Forman-Ricci curvature of vertices (Tables \ref{tab:ORFRedge}-\ref{tab:ORFRvertex} and Supplementary Tables S2-S3).

We have also compared the Ollivier-Ricci with Forman-Ricci and Augmented Forman-Ricci curvature of vertices in seventeen real-world networks. In several of the analyzed real-world networks, we find a moderate to high positive correlation between Ollivier-Ricci and Forman-Ricci curvature of vertices (Table \ref{tab:ORFRvertex} and Supplementary Table S3). Also, in most real-world networks analyzed here, we find a higher positive correlation between Augmented Forman-Ricci and Ollivier-Ricci curvature of vertices in comparison to Forman-Ricci and Ollivier-Ricci curvature of vertices (Table \ref{tab:ORFRvertex} and Supplementary Table S3). Thus, at the level of vertices, we observe a positive correlation between Ollivier-Ricci and Forman-Ricci curvature, especially, the augmented version, in many networks (Table \ref{tab:ORFRvertex} and Supplementary Table S2).

Importantly, we find that the correlation between Ollivier-Ricci and Forman-Ricci curvature of vertices is higher than Ollivier-Ricci and Forman-Ricci curvature of edges in most networks analyzed here (Tables \ref{tab:ORFRedge}-\ref{tab:ORFRvertex} and Supplementary Tables S2-S3). An intuitive explanation consists in the following observation. For the curvature of a vertex $v_0$ in an unweighted network, we average over all edges $(v_0,v)$ that have that vertex as one of its endpoints. Therefore, the Forman-Ricci curvature of each edge $(v_0,v)$ with vertex $v_0$ as one of its endpoint in an unweighted network has the form, $4-\deg v_0 -\deg v$ (see Eq. \ref{eq:Forman2}), and the Forman-Ricci curvature of all such edges $(v_0,v)$ share the term $\deg v_0$ which decreases the variance. For example, we even find a high positive correlation between Ollivier-Ricci and Forman-Ricci curvature of vertices in Email communication network where only a weak positive correlation exists between Ollivier-Ricci and Forman-Ricci curvature of edges (Tables \ref{tab:ORFRedge}-\ref{tab:ORFRvertex} and Supplementary Tables S2-S3). In a nut shell, although the two discretizations of Ricci curvature, Ollivier-Ricci and Forman-Ricci, capture different geometrical properties, our empirical analysis intriguingly finds a high positive correlation in many networks, especially, real-world networks. Deeper investigations in future are needed to better understand this empirically observed correlation between Ollivier-Ricci and Forman-Ricci curvature in many networks.


\begin{table}[ht]
\caption{Comparison of Ollivier-Ricci curvature (OR) with Forman-Ricci curvature (FR) or Augmented Forman-Ricci curvature (AFR) of edges in model and real networks. In this table, we list the Spearman correlation between the edge curvatures. In case of model networks, the reported correlation is mean (rounded off to two decimal places) over a sample of 100 networks generated with specific input parameters. Supplementary Table S2 also contains results from additional analysis of model networks with an expanded set of chosen input parameters. Moreover, Supplementary Table S2 also lists the Pearson correlation between the edge curvatures in model and real networks.}
\label{tab:ORFRedge}
\begin{tabular}{|l|c|c|}
\hline
{\textbf{\small Network}} & {\textbf{\small OR versus FR of edges}} & {\textbf{\small OR versus AFR of edges}} \\
\hline
\textbf{Model networks}                &            &   \\
\small{ER model with $n=1000$, $p=0.003$}   & 0.89  & 0.90  \\
\small{ER model with $n=1000$, $p=0.007$}   & 0.39  & 0.43  \\
\small{ER model with $n=1000$, $p=0.01$}   &  -0.03 & 0.04  \\
\small{WS model with $n=1000$, $k=2$ and $p=0.5$}   &  0.92 & 0.92  \\
\small{WS model with $n=1000$, $k=8$ and $p=0.5$}   &  0.18   & 0.70  \\
\small{WS model with $n=1000$, $k=10$ and $p=0.5$}   & 0.10   & 0.69  \\
\small{BA model with $n=1000$, $m=2$}   &  0.74      & 0.74  \\
\small{BA model with $n=1000$, $m=4$}   &  0.33      & 0.36  \\
\small{BA model with $n=1000$, $m=5$}   &  0.13      & 0.16  \\
\small{HGG model with $n=1000$, $k=3$, $\gamma=2$, $T=0$}   &  0.78      & 0.66  \\
\small{HGG model with $n=1000$, $k=5$, $\gamma=2$, $T=0$}   &  0.82      & 0.76  \\
\small{HGG model with $n=1000$, $k=10$, $\gamma=2$, $T=0$}   & 0.85       & 0.87  \\
\hline
\textbf{Real networks}                 &            &   \\
\small{Autonomous systems}             &    0.43    &   0.42    \\
\small{PGP}                            & 	0.32	&   0.83	\\
\small{US Power Grid}                  & 	0.60	&   0.76	\\
\small{Astrophysics co-authorship}     & 	0.25	&   0.70	\\
\small{Chicago Road}                   & 	0.98	&   0.98	\\
\small{Yeast protein interactions}     & 	0.70	&   0.74	\\
\small{Euro Road}                      & 	0.81	&   0.88	\\
\small{Human protein interactions}     & 	0.48	&   0.52	\\
\small{Hamsterster friendship}         & 	0.23	&   0.30	\\
\small{Email communication}            & 	0.19	&   0.53	\\
\small{PDZ domain interactions}        & 	0.72	&   0.71	\\
\small{Adjective-Noun adjacency}       & 	0.15	&   0.35	\\
\small{Dolphin}                        & 	0.07	&   0.71	\\
\small{Contiguous US States}           & 	0.68	&   0.91	\\
\small{Zachary karate club}            & 	0.75	&   0.81	\\
\small{Jazz musicians}                 & 	0.11	&   0.90	\\
\small{Zebra}                          & 	-0.04	&   0.62	\\
\hline
\end{tabular}
\end{table}

\begin{table}[ht]
\caption{Comparison of Ollivier-Ricci curvature (OR) with Forman-Ricci curvature (FR) or Augmented Forman-Ricci curvature (AFR) of vertices in model and real networks. In this table, we list the Spearman correlation between the vertex curvatures. In case of model networks, the reported correlation is mean (rounded off to two decimal places) over a sample of 100 networks generated with specific input parameters. Supplementary Table S3 also contains results from additional analysis of model networks with an expanded set of chosen input parameters. Moreover, Supplementary Table S3 also lists the Pearson correlation between the vertex curvatures in model and real networks.}
\label{tab:ORFRvertex}
\begin{tabular}{|l|c|c|}
\hline
{\textbf{\small Network}} & {\textbf{\small OR versus FR of vertices}} & {\textbf{\small OR versus AFR of vertices}} \\
\hline
\textbf{Model networks}                &            &   \\
\small{ER model with $n=1000$, $p=0.003$}   &  0.97     & 0.97  \\
\small{ER model with $n=1000$, $p=0.007$}   &  0.97     & 0.97  \\
\small{ER model with $n=1000$, $p=0.01$}   &   0.96    &  0.96 \\
\small{WS model with $n=1000$, $k=2$ and $p=0.5$}   & 0.90    & 0.90  \\
\small{WS model with $n=1000$, $k=8$ and $p=0.5$}   & 0.80    & 0.93  \\
\small{WS model with $n=1000$, $k=10$ and $p=0.5$}   & 0.77    & 0.92  \\
\small{BA model with $n=1000$, $m=2$}   &  0.61      & 0.61  \\
\small{BA model with $n=1000$, $m=4$}   &  0.59      & 0.60  \\
\small{BA model with $n=1000$, $m=5$}   &  0.63      & 0.64  \\
\small{HGG model with $n=1000$, $k=3$, $\gamma=2$, $T=0$}   &  0.48      & 0.57  \\
\small{HGG model with $n=1000$, $k=5$, $\gamma=2$, $T=0$}   &  0.34      & 0.41  \\
\small{HGG model with $n=1000$, $k=10$, $\gamma=2$, $T=0$}   & 0.09       & 0.13  \\
\hline
\textbf{Real networks}                 &            &   \\
\small{Autonomous systems}             &  0.64    &   0.64    \\
\small{PGP}                            & 	0.37	&   0.74	\\
\small{US Power Grid}                  & 	0.68	&   0.82	\\
\small{Astrophysics co-authorship}     & 	0.43	&   0.78	\\
\small{Chicago Road}                   & 	0.96	&   0.96	\\
\small{Yeast protein interactions}     & 	0.85	&   0.92	\\
\small{Euro Road}                      & 	0.90	&   0.92	\\
\small{Human protein interactions}     & 	0.83	&   0.84	\\
\small{Hamsterster friendship}         & 	0.85	&   0.86	\\
\small{Email communication}            & 	0.79	&   0.86	\\
\small{PDZ domain interactions}        & 	0.91	&   0.91	\\
\small{Adjective-Noun adjacency}       & 	0.47	&   0.50	\\
\small{Dolphin}                        & 	0.04	&   0.49	\\
\small{Contiguous US States}           & 	0.61	&   0.89	\\
\small{Zachary karate club}            & 	0.24	&   0.70	\\
\small{Jazz musicians}                 & 	-0.79	&   0.01	\\
\small{Zebra}                          & 	-0.72	&   0.99	\\
\hline
\end{tabular}
\end{table}

\subsection{Comparison of Forman-Ricci and Ollivier-Ricci curvature with other edge-based measures}

We emphasize that Ollivier-Ricci and Forman-Ricci curvature are edge-based measures of complex networks. We compared Ollivier-Ricci, Forman-Ricci and Augmented Forman-Ricci curvature with three other edge-based measures, edge betweenness centrality \cite{Freeman1977,Girvan2002,Newman2010}, embeddedness \cite{Marsden1984} and dispersion \cite{Backstrom2014}, for complex networks. Edge betweenness centrality \cite{Freeman1977,Girvan2002,Newman2010} measures the number of shortest paths that pass through an edge in a network. Edge betweenness centrality can be used to identify bottlenecks for flows in network. Embeddedness \cite{Marsden1984} of an edge quantifies the number of neighbors that are shared by the two vertices anchoring the edge under consideration in the network. Embeddedness is a measure to quantify the strength of ties in social networks \cite{Marsden1984}. Dispersion \cite{Backstrom2014} quantifies the extent to which the neighbours of the two vertices anchoring an edge are not themselves well connected. Dispersion is a measure to predict romantic relationships in social networks \cite{Backstrom2014}.

In model networks, we find that Ollivier-Ricci, Forman-Ricci and Augmented Forman-Ricci curvature have significant negative correlation with edge betweenness centrality (Table \ref{tab:edge} and Supplementary Table S4). In most real networks considered here, we find that Ollivier-Ricci curvature has moderate to high negative correlation with edge betweenness centrality while Forman-Ricci curvature has a weak to moderate negative correlation with edge betweenness centrality (Table \ref{tab:edge} and Supplementary Table S4). Moreover, in most real networks considered here, we observe a higher negative correlation between Ollivier-Ricci curvature and edge betweenness centrality in comparison to Forman-Ricci curvature and edge betweenness centrality (Table \ref{tab:edge} and Supplementary Table S4). This may be explained by the fact that Ollivier-Ricci curvature is also affected by cycles of length 3, 4 and 5 containing the two vertices of an edge, and these are relevant for edge betweenness centrality. Interestingly, in real networks considered here, the Augmented Forman-Ricci curvature in comparison to Forman-Ricci curvature has much higher negative correlation with edge betweenness centrality (Table \ref{tab:edge} and Supplementary Table S4). Our results suggest that the augmented version of Forman-Ricci curvature which also accounts for two-dimensional simplicial complexes arising in graphs is better suited for analysis of complex networks.

In both model and real networks considered here, we find no consistent relationship between Ollivier-Ricci, Forman-Ricci, or Augmented Forman-Ricci curvature of an edge and embeddedness (Table \ref{tab:edge} and Supplementary Table S4). Similarly, In both model and real networks considered here, we find no consistent relationship between Ollivier-Ricci, Forman-Ricci, or Augmented Forman-Ricci curvature of an edge and dispersion (Table \ref{tab:edge} and Supplementary Table S4). In summary, the two discrete notions of Ricci curvatures are negatively correlated to edge betweeness centrality but have no consistent relationship with embeddedness or dispersion in analyzed networks.

\begin{table}[ht]
\caption{ Comparison of Ollivier-Ricci curvature (OR), Forman-Ricci curvature (FR) and Augmented Forman-Ricci curvature (AFR) of edges with other edge-based measures, edge betweenness centrality (EBC), embeddedness (EMB) and dispersion (DIS), in model and real networks. In this table, we list the Spearman correlation between the edge-based measures. In case of model networks, the reported correlation is mean (rounded off to two decimal places) over a sample of 100 networks generated with specific input parameters. Supplementary Table S4 also contains results from additional analysis of model networks with an expanded set of chosen input parameters. Moreover, Supplementary Table S4 also lists the Pearson correlation between the edge-based measures in model and real networks.}
\label{tab:edge}
\begin{tabular}{|l|c|c|c|c|c|c|c|c|c|}
\hline
{\textbf{\small Network}} & \multicolumn{3}{c|}{\textbf{\small OR versus}} & \multicolumn{3}{c|}{\textbf{\small FR versus}} & \multicolumn{3}{c|}{\textbf{\small AFR versus}} \\
\cline{2-10}
& \textbf{\small EBC} & \textbf{\small EMB} & \textbf{\small DIS}  & \textbf{\small EBC} & \textbf{\small EMB} & \textbf{\small DIS} & \textbf{\small EBC} & \textbf{\small EMB} & \textbf{\small DIS}\\
\hline
\textbf{Model networks}                     &   &   &   &   &   &   &   &   &   \\
\small{ER model with $n=1000$, $p=0.003$}  &  -0.86 & 0.08  & 0.00  & -0.81  & -0.07  & 0.00  & -0.82  & 0.04  & 0.00  \\
\small{ER model with $n=1000$, $p=0.007$}   & -0.53  & 0.25  & 0.05  & -0.80  & -0.11  & -0.03  & -0.82  & 0.06  & 0.02  \\
\small{ER model with $n=1000$, $p=0.01$}   & -0.34  & 0.32  & 0.10  & -0.76  & -0.13  & -0.05  & -0.79  & 0.07  & 0.03  \\
\small{WS model with $n=1000$, $k=2$ and $p=0.5$}   & -0.75  & 0.00  & 0.00  & -0.57  & 0.00  & 0.00  & -0.57  & 0.00  & 0.00  \\
\small{WS model with $n=1000$, $k=8$ and $p=0.5$}   & -0.85  & 0.79  & 0.44  & -0.52  & -0.05  & -0.08  & -0.89  & 0.68  & 0.42  \\
\small{WS model with $n=1000$, $k=10$ and $p=0.5$}   & -0.87  & 0.82  & 0.49  & -0.45  & -0.05  & -0.07  & -0.89  & 0.73  & 0.47  \\
\small{BA model with $n=1000$, $m=2$}   & -0.73  & -0.09  & -0.11  & -0.76  & -0.30  & -0.16  & -0.77  & -0.26  & -0.15  \\
\small{BA model with $n=1000$, $m=4$}   & -0.45  & 0.18  & 0.14  & -0.83  & -0.48  & -0.35  & -0.84  & -0.43  & -0.33  \\
\small{BA model with $n=1000$, $m=5$}   & -0.30  & 0.30  & 0.25  & -0.85  & -0.54  & -0.41  & -0.86  & -0.48  & -0.39  \\
\small{HGG model with $n=1000$, $k=3$, $\gamma=2$, $T=0$}   & -0.47  & -0.30  & -0.15  & -0.67  & -0.04  & -0.18  & -0.76  & 0.27  & -0.07  \\
\small{HGG model with $n=1000$, $k=5$, $\gamma=2$, $T=0$}   &  -0.62 & -0.20  & -0.13  & -0.73  & -0.08  & -0.17  & -0.81  & 0.20  & -0.10  \\
\small{HGG model with $n=1000$, $k=10$, $\gamma=2$, $T=0$}   & -0.78  & -0.03  & -0.06  & -0.79  & -0.15  & -0.12  & -0.87  & 0.14  & -0.08  \\
\hline
\textbf{Real networks}                 &   &   &   &   &   &   &   &   &   \\
\small{Autonomous systems}             &  -0.17 & -0.37 & -0.25 & -0.26 & -0.44  & -0.18  & -0.27  &  -0.41 & -0.16  \\
\small{PGP}                            &  -0.64 & 0.20  & -0.13 & 0.11  & -0.69  & -0.17  & -0.56  &  0.21  & -0.15  \\
\small{US Power Grid}                  &  -0.61 & 0.16  & 0.06  & -0.26 & -0.41  & -0.19  & -0.45  &  0.09  & 0.04  \\
\small{Astrophysics co-authorship}     &  -0.78 & 0.47  & -0.16 & -0.23 & -0.58  & -0.23  & -0.63  &  0.07  & -0.27  \\
\small{Chicago Road}                   &  -0.65 & 0.00  & 0.00  & -0.65 & 0.00   & 0.00   & -0.65  &  0.00  & 0.00  \\
\small{Yeast protein interactions}     &  -0.83 & 0.06  & -0.01 & -0.52 & -0.15  & -0.13  & -0.59  &  0.14  & 0.00  \\
\small{Euro Road}                      &  -0.54 & 0.05  & 0.02  & -0.40 & -0.31  & -0.07  & -0.43  &  0.00  & 0.03  \\
\small{Human protein interactions}     &  -0.46 & 0.07  & 0.01  & -0.38 & -0.22  & -0.19  & -0.43  &  -0.07 & -0.10  \\
\small{Hamsterster friendship}         &  -0.53 & 0.12  & 0.00  & -0.35 & -0.61  & -0.40  & -0.42  &  -0.47 & -0.32  \\
\small{Email communication}            &  -0.61 & 0.55  & 0.24  & -0.32 & -0.45  & -0.41  & -0.57  &  0.01  & -0.16  \\
\small{PDZ domain interactions}        &  -0.79 & -0.04 & 0.00  & -0.55 & -0.02  & 0.00   & -0.55  &  0.06  & 0.00  \\
\small{Adjective-Noun adjacency}       &  -0.51 & 0.22  & 0.09  & -0.42 & -0.72  & -0.55  & -0.57  &  -0.42 & -0.37  \\
\small{Dolphin}                        &  -0.66 & 0.51  & 0.28  & 0.11  & -0.58  & -0.21  & -0.61  &  0.59  & 0.31  \\
\small{Contiguous US States}           &  -0.68 & -0.10 & -0.15 & -0.49 & -0.72  & -0.71  & -0.64  &  -0.03 & -0.08  \\
\small{Zachary karate club}            &  -0.79 & 0.10  & -0.06 & -0.64 & -0.29  & -0.37  & -0.80  &  0.43  & 0.14  \\
\small{Jazz musicians}                 &  -0.84 & 0.57  & -0.03 & -0.22 & -0.66  & -0.18  & -0.76  &  0.47  & -0.05  \\
\small{Zebra}                          &  -0.94 & 0.52  & 0.13  & 0.04  & -0.71  & -0.15  & -0.65  &  0.97  & 0.09  \\
\hline
\end{tabular}
\end{table}

\begin{table}[ht]
\caption{ Comparison of Ollivier-Ricci curvature (OR), Forman-Ricci curvature (FR) and Augmented Forman-Ricci curvature (AFR) of vertices with other vertex-based measures, degree, betweenness centrality (BC) and clustering coefficient (CC), in model and real networks. In this table, we list the Spearman correlation between the vertex-based measures. In case of model networks, the reported correlation is mean (rounded off to two decimal places) over a sample of 100 networks generated with specific input parameters. Supplementary Table S5 also contains results from additional analysis of model networks with an expanded set of chosen input parameters. Moreover, Supplementary Table S5 also lists the Pearson correlation between the vertex-based measures in model and real networks.}
\label{tab:vertex}
\begin{tabular}{|l|c|c|c|c|c|c|c|c|c|}
\hline
{\textbf{\small Network}} & \multicolumn{3}{c|}{\textbf{\small OR versus}} & \multicolumn{3}{c|}{\textbf{\small FR versus}} & \multicolumn{3}{c|}{\textbf{\small AFR versus}} \\
\cline{2-10}
& \textbf{\small Degree} & \textbf{\small BC}  & \textbf{\small CC} & \textbf{\small Degree} & \textbf{\small BC}  & \textbf{\small CC} & \textbf{\small Degree} & \textbf{\small BC}  & \textbf{\small CC} \\
\hline
\textbf{Model networks}                     &   &   &   &   &   &   &   &   &   \\
\small{ER model with $n=1000$, $p=0.003$}  &  -0.94 & -0.94  & -0.07  & -0.94  & -0.94  & -0.13  & -0.94  & -0.94  & -0.08  \\
\small{ER model with $n=1000$, $p=0.007$}   & -0.98  & -0.98  & -0.18  & -0.99  & -0.98  & -0.26  & -0.99  & -0.98  & -0.21  \\
\small{ER model with $n=1000$, $p=0.01$}   & -0.98  & -0.98  & -0.16  & -0.99  & -0.98  & -0.25  & -0.99  & -0.98  & -0.21  \\
\small{WS model with $n=1000$, $k=2$ and $p=0.5$}   & -0.71  & -0.82  & 0.00  & -0.75  & -0.73  & 0.00  & -0.75  & -0.73  & 0.00  \\
\small{WS model with $n=1000$, $k=8$ and $p=0.5$}   & -0.81  & -0.96  & 0.51  & -0.98  &  -0.91 & 0.05  & -0.91  & -0.98  & 0.38  \\
\small{WS model with $n=1000$, $k=10$ and $p=0.5$}   & -0.79  & -0.95  & 0.57  & -0.99  & -0.91  & 0.09  & -0.92  & -0.98  & 0.41  \\
\small{BA model with $n=1000$, $m=2$}   & -0.90  & -0.90  & -0.18  & -0.59  & -0.77  & -0.39  & -0.59  & -0.78  & -0.37  \\
\small{BA model with $n=1000$, $m=4$}   & -0.94  & -0.88  & -0.08  & -0.73  & -0.84  &  -0.49 &  -0.73 & -0.85  & -0.45  \\
\small{BA model with $n=1000$, $m=5$}   & -0.94  & -0.90  & -0.05  & -0.78  & -0.85  & -0.40  & -0.79  & -0.86  & -0.37  \\
\small{HGG model with $n=1000$, $k=3$, $\gamma=2$, $T=0$}   & -0.28  & -0.30  & -0.14  & -0.86  & -0.60  & -0.45  & -0.79  & -0.58  & -0.37  \\
\small{HGG model with $n=1000$, $k=5$, $\gamma=2$, $T=0$}   & -0.15  & -0.17  & -0.03  & -0.89  & -0.61  & -0.21  & -0.85  & -0.60  & -0.18  \\
\small{HGG model with $n=1000$, $k=10$, $\gamma=2$, $T=0$}   & 0.06  & -0.06  & 0.01  & -0.93  & -0.68  & 0.31  & -0.91  & -0.66  & 0.30  \\
\hline
\textbf{Real networks}                 &   &   &   &   &   &   &   &   &   \\
\small{Autonomous systems}             & -0.85  & -0.70  & -0.39 & -0.51  & -0.38  & -0.55 & -0.50  & -0.38  & -0.55  \\
\small{PGP}                            & -0.12  & -0.49  & 0.29  & -0.73  & -0.51  & -0.51 & -0.35  & -0.46  & -0.05  \\
\small{US Power Grid}                  & -0.68  & -0.80  & 0.03  & -0.79  & -0.62  & -0.49 & -0.69  & -0.68  & -0.13  \\
\small{Astrophysics co-authorship}     & -0.39  & -0.72  & 0.62  & -0.95  & -0.64  & 0.25  & -0.64  & -0.66  & 0.41  \\
\small{Chicago Road}                   & -0.33  & -0.34  & 0.00  & -0.42  & -0.42  & 0.00  & -0.42  & -0.42  & 0.00  \\
\small{Yeast protein interactions}     & -0.54  & -0.67  & -0.05 & -0.57  & -0.56  & -0.33 & -0.45  & -0.54  & -0.07  \\
\small{Euro Road}                      & -0.82  & -0.75  & -0.22 & -0.82  & -0.64  & -0.38 & -0.80  & -0.65  & -0.24  \\
\small{Human protein interactions}     & -0.77  & -0.78  & -0.23 & -0.71  & -0.65  & -0.43 & -0.67  & -0.64  & -0.34  \\
\small{Hamsterster friendship}         & -0.87  & -0.87  & -0.30 & -0.92  & -0.76  & -0.45 & -0.91  & -0.76  & -0.42  \\
\small{Email communication}            & -0.80  & -0.88  & 0.06  & -0.97  & -0.87  & -0.31 & -0.93  & -0.88  & -0.19  \\
\small{PDZ domain interactions}        & -0.50  & -0.58  & -0.12 & -0.62  & -0.64  & -0.14 & -0.61  & -0.64  & -0.09  \\
\small{Adjective-Noun adjacency}       & -0.57  & -0.76  & 0.07  & -0.96  & -0.84  & -0.50 & -0.95  & -0.84  & -0.45  \\
\small{Dolphin}                        & -0.04  & -0.39  & 0.44  & -0.98  & -0.77  & -0.45 & -0.73  & -0.72  & -0.04  \\
\small{Contiguous US States}           & -0.59  & -0.74  & 0.71  & -0.98  & -0.82  & 0.55  & -0.78  & -0.79  & 0.70  \\
\small{Zachary karate club}            & 0.10   & -0.09  & 0.35  & -0.84  & -0.76  & 0.40  & -0.47  & -0.60  & 0.52  \\
\small{Jazz musicians}                 & 0.78   & 0.34   & 0.08  & -0.99  & -0.72  & 0.33  & -0.49  & -0.56  & 0.56  \\
\small{Zebra}                          & 0.78   & 0.35   & -0.33 & -0.94  & -0.73  & 0.70  & 0.76   & 0.33   & -0.31  \\
\hline
\end{tabular}
\end{table}
\subsection{Comparison of Forman-Ricci and Ollivier-Ricci curvature with vertex-based measures}

We compared Ollivier-Ricci, Forman-Ricci and Augmented Forman-Ricci curvature of vertices with three other vertex-based measures, degree, betweenness centrality \cite{Freeman1977,Newman2010} and clustering coefficient \cite{Holland1971,Watts1998}, in a network. Vertex degree gives the number of edges incident to that vertex in a network. Betweennness centrality \cite{Freeman1977,Newman2010} of a vertex quantifies the fraction of shortest paths between all pairs of vertices in the network that pass through that vertex. The clustering coefficient \cite{Holland1971,Watts1998} of a vertex quantifies the number of edges that are realized between the neighbours of the vertex divided by the number of edges that could possibly exist between the neighbours of the vertex in the network. We remark that the clustering coefficient has been proposed as a measure to quantify the curvature of networks \cite{Eckmann2002}.

Not surprisingly, we find that Ollivier-Ricci, Forman-Ricci or Augmented Forman-Ricci curvature of vertices have high negative correlation with degree in most model as well as real networks analyzed here (Table \ref{tab:vertex} and Supplementary Table S5). After all, the vertex degree is intrinsic in the definition of the Ollivier-Ricci or Forman-Ricci curvature of a vertex as it appears implicitly in the sum over adjacent edges in the defining formula. Similarly, in model as well as real networks analyzed here, we find that Ollivier-Ricci, Forman-Ricci or Augmented Forman-Ricci curvature of vertices have high negative correlation with betweenness centrality (Table \ref{tab:vertex} and Supplementary Table S5). In contrast, we do not find any consistent relationship between Ollivier-Ricci, Forman-Ricci or Augmented Forman-Ricci curvature of vertices and clustering coefficient in model and real networks analyzed here (Table \ref{tab:vertex} and Supplementary Table S5).

\begin{figure}
\includegraphics[width=.7\columnwidth]{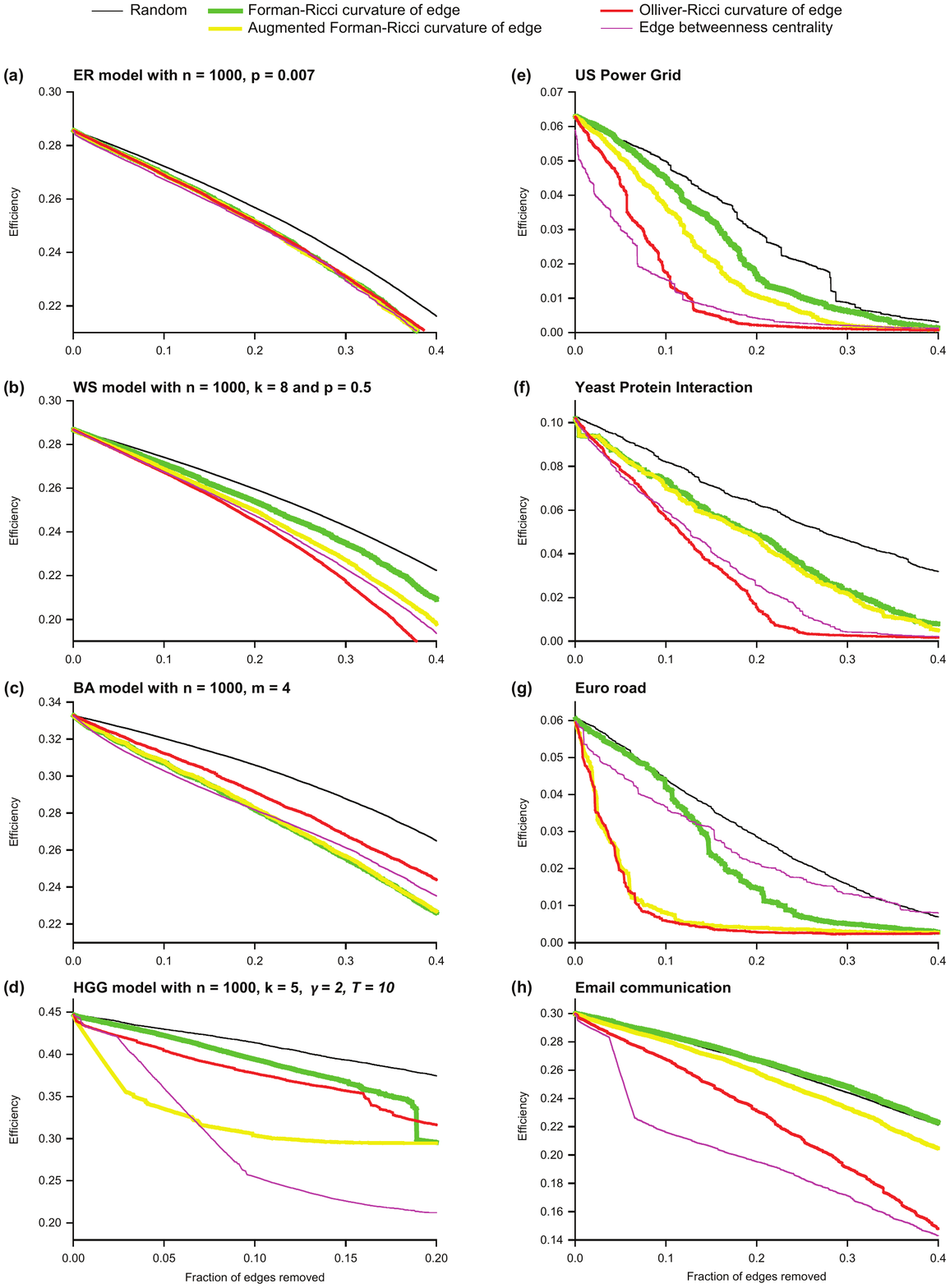}
\caption{Communication efficiency as a function of the fraction of edges removed in model and real networks. (a) Erd\"{o}s-R\`{e}nyi (ER) model. (b) Watts-Strogratz (WS) model. (c) Barab\`{a}si-Albert (BA) model. (d) Hyberbolic random geometric graph (HGG) model. (e) US Power Grid. (f) Yeast protein interactions. (g) Euro road. (h) Email communication. }
\label{fig:rob_edge}
\end{figure}

\begin{figure}
\includegraphics[width=.7\columnwidth]{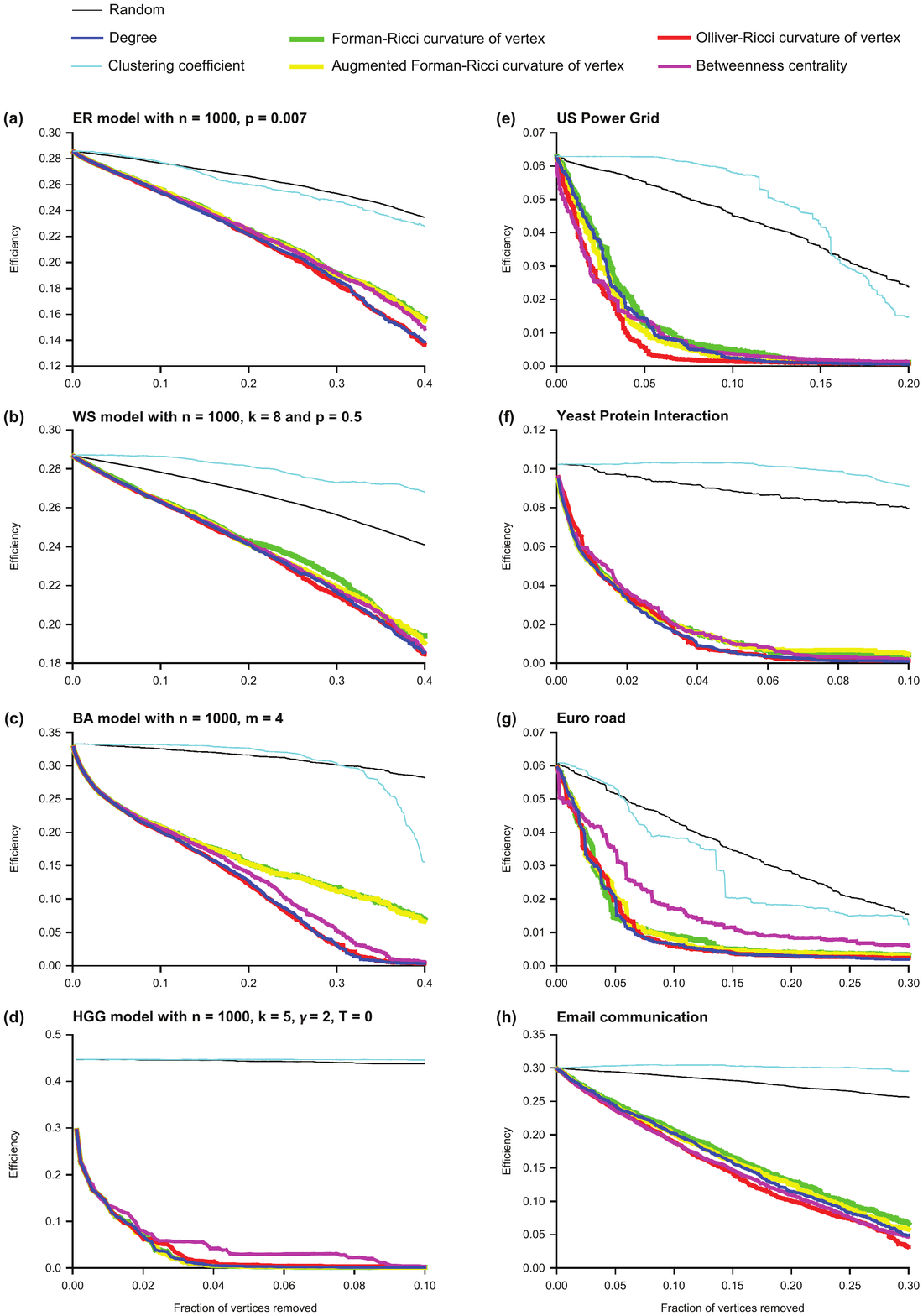}
\caption{Communication efficiency as a function of the fraction of vertices removed in model and real networks. (a) Erd\"{o}s-R\`{e}nyi (ER) model. (b) Watts-Strogratz (WS) model. (c) Barab\`{a}si-Albert (BA) model. (d) Hyberbolic random geometric graph (HGG) model. (e) US Power Grid. (f) Yeast protein interactions. (g) Euro road. (h) Email communication.}
\label{fig:rob_node}
\end{figure}
\subsection{Relative importance of Forman-Ricci and Ollivier-Ricci curvature for topological robustness of networks}

We employ a global network measure, communication efficiency \cite{Latora2001}, to quantify the effect of removing edges or vertices on the large-scale connectivity of networks. Communication efficiency $E$ of a graph $G$ is given by:
\begin{equation}
E = \frac{1}{n(n-1)}\sum_{i < j \in V(G)}\frac{1}{d_{ij}},
\end{equation}
where $d_{ij}$ denotes the shortest path between the pair of vertices $i$ and $j$, $n$ is the number of vertices in the graph, and $V(G)$ denotes the set of vertices in the graph. Note that communication efficiency captures the resilience of a network to failure in the face of perturbations, as it essentially identifies locally with the clustering coefficient and globally with the inverse of the characteristic path length.

We investigated the relative importance of Ollivier-Ricci, Forman-Ricci or Augmented Forman-Ricci curvature of edges for the large-scale connectivity of networks by removing edges based on the following criteria: random order, increasing order of the Forman-Ricci curvature of an edge, increasing order of the Augmented Forman-Ricci curvature of an edge, increasing order of the Ollivier-Ricci curvature of an edge, and decreasing order of edge betweenness centrality. In both model and real networks, we find that removing edges based on increasing order of Ollivier-Ricci curvature or increasing order of Forman-Ricci curvature or increasing order of Augmented Forman-Ricci curvature or decreasing order of edge betweenness centrality leads to faster disintegration in comparison to the random removal of edges (Figure \ref{fig:rob_edge}). Furthermore, in most cases, removing edges based on increasing order of Ollivier-Ricci curvature or decreasing order of edge betweenness centrality typically leads to faster disintegration in comparison to removing edges based on increasing order of Forman-Ricci curvature (Figure \ref{fig:rob_edge}). We remark that both Ollivier-Ricci curvature of an edge and edge betweenness centrality are global measures while Forman-Ricci curvature of an edge is a local measure dependent on nearest neighbors of an edge.

We also investigated the relative importance of Ollivier-Ricci, Forman-Ricci or Augmented Forman-Ricci curvature of vertices for the large-scale connectivity of networks by removing vertices based on the following criteria: random order, increasing order of the Forman-Ricci curvature of a vertex, increasing order of the Augmented Forman-Ricci curvature of a vertex, increasing order of the Ollivier-Ricci curvature of a vertex, decreasing order of betweenness centrality of a vertex, decreasing order of vertex degree, and decreasing order of clustering coefficient of a vertex. In both model and real networks, we find that removing vertices based on increasing order of Ollivier-Ricci curvature or increasing order of Forman-Ricci curvature or increasing order of Augmented Forman-Ricci curvature or decreasing order of betweenness centrality or decreasing order of degree leads to faster disintegration in comparison to the random removal of vertices (Figure \ref{fig:rob_node}). Furthermore, in most model as well as real networks, removing vertices based on increasing order of Ollivier-Ricci curvature typically leads to faster disintegration in comparison to removing edges based on increasing order of Forman-Ricci curvature or on increasing order of Augmented Forman-Ricci curvature (Figure \ref{fig:rob_node}). Also, in most model as well as real networks, removing edges based on increasing order of Ollivier-Ricci curvature typically leads to at least slightly faster disintegration in comparison to removing edges based on any other measure (Figure \ref{fig:rob_node}). In summary, vertices or edges with highly negative Ollivier-Ricci curvature are found to be more important than vertices or edges with highly negative Forman-Ricci curvature for maintaining the large-scale connectivity of most networks analyzed here.


\section{Conclusions}

We have performed an empirical investigation of two discretizations of Ricci curvature, Ollivier's Ricci curvature and Forman's Ricci curvature, in a number of model and real-world networks. The two discretizations of Ricci curvature were derived using different theoretical considerations and methods, and thus, convey insights into quite different geometrical properties and behaviors of complex networks. Specifically, Ollivier-Ricci curvature captures clustering and coherence in networks while Forman-Ricci curvature captures dispersal and topology. Moreover, in the context of weighted networks, Ollivier-Ricci curvature implicitly, by its very definition, relates to edge weights as probabilities, while Forman's Ricci curvature fundamentally views edge weights as abstractions of lengths, and vertex weights as, for instance, concentrated area measures. This suggests that Ollivier-Ricci curvature is intrinsically better suited to study probabilistic phenomenon on networks while Forman-Ricci curvature is better suited to investigate networks where edge weights correspond to distances. Still, our results obtained in a wide-range of both model and real-world networks, consistently demonstrate that the two types of Ricci curvature in many networks are highly correlated. The immediate benefit of this realization is that one can compute Forman-Ricci curvature in large networks to gain some first insight into the computationally much more demanding Ollivier-Ricci curvature. Furthermore, the state of the art computational implementation of the Ollivier-Ricci curvature can handle only weights on edges rather than vertices in weighted networks. In addition, while computing the Ollivier-Ricci curvature of an edge in a weighted network, a necessary step is the normalization of the neighboring edge weights. In contrast, the mathematical definition of the Forman-Ricci curvature can incorporate any set of positive weights, placed simultaneously at the vertices and the edges. Furthermore, the Augmented Forman-Ricci curvature can also account for higher-dimensional simplicial complexes, thus making it a natural and simple to employ tool for the understanding networks with explicit geometric structure, especially, hyper-networks. Therefore, our empirical observations on the correlation between these two different notions of Ricci curvature in networks warrant deeper investigation in the future.

We remark that while the present manuscript was under final stages of submission, a preprint \cite{Pouryahya2017} devoted to comparison problem in biological networks appeared on Arxiv server, independently from our present study.


\section*{Acknowledgments}

We thank the anonymous reviewers for their constructive comments which have helped improve the manuscript. E.S. and A.S. thank the Max Planck Institute for Mathematics in the Sciences, Leipzig, for their warm hospitality. A.S. would like to acknowledge support from Max Planck Society, Germany, through the award of a Max Planck Partner Group in Mathematical Biology.


%

\end{document}